\documentclass[10pt]{amsart}

\usepackage{amsmath,amssymb,amsfonts,amsthm,graphics, verbatim,psfrag}
\usepackage{latexsym,epsfig,verbatim}
\usepackage{textcomp}
\usepackage{float}
\usepackage{lineno}

\usepackage[all,2cell,dvips]{xy} 

\include{diagram}

\newcommand{\nc}{\newcommand}
\nc{\ens}{\ensuremath}
\nc{\dmo}{\DeclareMathOperator}
\nc{\nt}{\newtheorem}

\nt{thm}{Theorem}
\nt{mthm}{Main Theorem}
\nt{cor}[thm]{Corollary}
\nt{minithm}[thm]{Proposition}
\nt{lem}[thm]{Lemma}
\nt{step}[thm]{Step}
\nt{fact}[thm]{Fact}
\nt{q}{Question}
\nt{rem}{Remark}

\nc{\prop}{Proposition}

\dmo{\rk}{rk}
\nc{\cent}{\ens{Z}}
\nc{\K}{\mathcal K}
\nc{\I}{\mathcal I}
\nc{\C}{\mathcal C}
\nc{\Csep}{\C_s}
\nc{\ps}{\phi_\star}
\nc{\pair}{{\mathcal P}}
\nc{\curv}{{\mathcal C}}
\dmo{\Comm}{Comm}
\nc{\commk}{\comm(\K)}
\nc{\fin}{G}

\nc{\bc}{\mathcal{B}}

\nc{\pst}{\widehat{\phi}_{\star}}

\nc{\uc}{\al\text{-}\be\text{-}\ga}
\nc{\ucp}{\uc'}
\nc{\ucpp}{\al'\text{-}\be\text{-}\ga'}
\nc{\ucppp}{\al'\text{-}\be'\text{-}\ga'}
\nc{\newc}{$\al'$-$\be$-$\ga'$ \hskip .01in}
\nc{\tri}{$\be^-\al\de$ \hskip .01in}
\nc{\bad}{$\be^\pm\al\de$ \hskip .01in}
\nc{\bae}{$\be^\pm\al\ep$ \hskip .01in}
\nc{\bcd}{$\be^\pm\ga\de$ \hskip .01in}
\nc{\bce}{$\be^\pm\ga\ep$ \hskip .01in}
\nc{\adce}{$\al\de\ga\ep$ \hskip .01in}
\nc{\bdce}{$\be^\pm\de\ga\ep$ \hskip .01in}
\nc{\bcea}{$\be^\pm\ga\ep\al$ \hskip .01in}
\nc{\bead}{$\be^\pm\ep\al\de$ \hskip .01in}
\nc{\badc}{$\be^\pm\al\de\ga$ \hskip .01in}
\nc{\badce}{$\be^\pm\al\de\ga\ep$ \hskip .01in}
\nc{\bcead}{$\be^\pm\ga\al\ep\de$ \hskip .01in}
\nc{\bcdae}{$\be^\pm\ga\de\al\ep$ \hskip .01in}
\nc{\baecd}{$\be^\pm\al\ep\ga\de$ \hskip .01in}
\nc{\bottom}{$\be^-$}
\nc{\topp}{$\be^+$}

\nc{\td}{\text{-}}

\nc{\al}{\ens{\alpha}}
\nc{\be}{\ens{\beta}}
\nc{\ga}{\ens{\gamma}}
\nc{\de}{\ens{\delta}}
\nc{\ep}{\ens{\epsilon}}
\nc{\Ga}{\ens{\Gamma}}
\nc{\ta}{\ens{\tilde{a}}}
\nc{\tb}{\ens{\tilde{b}}}
\nc{\De}{\ens{\Delta}}

\dmo{\Mod}{Mod} 
\dmo{\pmcg}{PMod} 
\nc{\mods}{\ens{\Mod(S_g)}}
\dmo{\diff}{Diff}

\dmo{\homeop}{\ens{Homeo^+}}
\dmo{\homeopm}{\ens{Homeo^\pm}}
\dmo{\homeo}{\ens{Homeo}}

\dmo{\gin}{i} 

\dmo{\cc}{C}
\nc{\ccs}{\cc(S)}
\nc{\ccos}{\cc^0(S)}
\nc{\csep}{\cc_s(S)}

\nc{\csepo}{\cc_s^0(S)}

\nc{\tg}{\mathcal T(S_g)}
\nc{\T}{\mathcal T}
\nc{\tgo}{\mathcal T^{0
}(S)}
\nc{\tgeom}{\mathcal TG(S)}

\dmo{\si}{SI}

\dmo{\out}{Out}
\dmo{\Aut}{Aut}

\nc{\zt}{\ens{{\mathbb Z}_2}}

\nc{\pslz}{\ensuremath{PSL_2(\z)}}
\nc{\spz}{\ens{Sp(2g,\z)}}

\nc{\Z}{\ens{{\mathbb Z}}}
\nc{\br}{\ensuremath{{\mathbb R}}}

\dmo{\isom}{Isom} 

\nc{\lra}{\ens{\longrightarrow}}
\nc{\ra}{\ens{\rightarrow}}
\nc{\set}[1]{\{#1\}}
\nc{\genby}[1]{\langle #1 \rangle}

\nc{\bs}{\bigskip}
\nc{\bpf}{\begin{proof}}
\nc{\epf}{\end{proof}}
\nc{\p}[1]{{\bf #1} }

\nc{\bcom}{\begin{comment}}
\nc{\ecom}{\end{comment}}

\nc{\bl}{ \begin{list}{$\cdot$}{
\setlength{\leftmargin}{.25in}
\setlength{\rightmargin}{.5in}
\setlength{\parsep}{0.5ex plus .2ex minus 0ex}
\setlength{\itemsep}{0.2ex plus 0.2ex minus 0ex} 
}
}

\nc{\el}{\end{list}}

\nc{\pic}[2]{\begin{figure}[htb] \center{\leavevmode 
\epsfbox{#1.eps}} \caption{#2} \label{#1pic}\end{figure}} 

\nc{\pics}[3]{\epsfysize=#3 cm \begin{figure}[htb] 
\center{\leavevmode \epsfbox{#1.eps}} \caption{#2}  
\label{#1pic}\end{figure}}

\begin{document}

\input{epsf.sty}

\title{Addendum to: Commensurations of the Johnson Kernel}

\author{Tara E. Brendle}

\author{Dan Margalit} 

\address{Tara E. Brendle, Dept. of Mathematics, Louisiana State University, Baton Rouge, LA 70803-4918}

\email{brendle@math.lsu.edu}

\address{Dan Margalit, Dept. of Mathematics, University of Utah, 155 S. 1400 East, Salt Lake City, UT 84112}

\email{margalit@math.utah.edu}

\thanks{The first author was supported in part by NSF grant
  DMS-0606882 and the LSU Council on Research Summer Stipend Program.
  The second author was partially supported by NSF grant DMS-0707279
  and the NSF VIGRE program.}



\maketitle

\begin{center}\today\end{center}

\begin{abstract}Let $\K(S)$ be the subgroup of the extended mapping
class group, $\Mod(S)$, generated by Dehn twists about separating
curves.  In our earlier paper, we showed that $\Comm(\K(S)) \cong
\Aut(\K(S)) \cong \Mod(S)$ when $S$ is a closed, connected, orientable
surface of genus $g \geq 4$.  By modifying our original proof, we show
that the same result holds for $g \geq 3$, thus confirming Farb's
conjecture in all cases (the statement is not true for $g \leq 2$).
\end{abstract}

%%%
%%%
%%%

The purpose of this note is to extend the results of our paper
``Commensurations of the Johnson kernel'' to the lone remaining case.
We briefly review the notation and basic ideas before explaining the
improvement.  We refer to reader to that paper for further details
\cite{bm}.

Let $S=S_g$ denote a closed, connected, orientable, surface of genus $g$,
and let $\Mod(S)$ denote the extended mapping class group (orientation
reversing elements are allowed).  The Torelli group $\I(S)$ is the
subgroup of $\Mod(S)$ consisting of elements which act trivially on
$H_1(S,\Z)$, and the Johnson kernel $\K(S)$ is the subgroup of $\I(S)$ generated by Dehn twists about separating curves.

The {\em abstract commensurator} of a group $\Ga$, denoted $\Comm(\Ga)$, is the group of isomorphisms of finite index subgroups of $\Ga$ (under composition), with two such isomorphisms equivalent if they agree on a finite index subgroup of $\Ga$.  The product of $\phi : G \to H$ with $\psi : G' \to H'$ is a map defined on $\phi^{-1}(H \cap G')$.

We have the following theorem, which confirms a conjecture of Farb in
all cases \cite{bf}.  In the original paper, we stated and proved the
theorem for $g \geq 4$.

\begin{thm}\label{comm}
Let $g \geq 3$, and let $G$ be either $\I(S_g)$ or $\K(S_g)$.  We have 
\[ \Comm(G) \cong \Aut(G) \cong \Mod(S_g). \]
\end{thm}
For $g \geq 5$ and $G=\I(S_g)$, Theorem~\ref{comm} is due to
Farb--Ivanov \cite{fi}.  McCarthy--Vautaw proved that $\Aut(\I(S_g))
\cong \Mod(S_g)$ for $g\geq 3$ \cite{mv}.  Mess proved that $\I(S_2)=\K(S_2)$ 
is an infinitely generated free group,
so Theorem~\ref{comm} certainly does not hold in this case \cite{gm}.  Also, it
is a theorem of Dehn that $\I(S_1) = 1$.

Theorem~\ref{comm} is a consequence of the following more general theorem.

\begin{thm}\label{intorelli}
Let $g \geq 3$, let $H$ be a finite index subgroup of either $\I(S_g)$
or $\K(S_g)$.  Any injective homomorphism $\phi: H \to \I(S_g)$ is
induced by an element $f$ of $\Mod(S_g)$ in the sense that $\phi(h) =
fhf^{-1}$ for all $h \in H$.
\end{thm}

Theorem~\ref{intorelli} has various corollaries.  In particular, it
follows that finite index subgroups of $\I(S)$ and $\K(S)$ are
co-Hopfian, and that finite index subgroups of $\I(S)$
and $\K(S)$ are characteristic in $\I(S)$ up to conjugacy.

Our basic method, following Ivanov, is to translate
Theorem~\ref{intorelli} into a question about curve complexes.  The
complex of curves $\C(S)$ is the complex with vertices for isotopy
classes of simple closed curves in $S$ and simplices for disjointness.  The
complex of separating curves $\Csep(S)$ is the subcomplex spanned by
the separating curves.  Finally, the Torelli complex $\T(S)$  has vertices for
isotopy classes of separating curves and isotopy classes of bounding
pairs in $S$, and simplices for disjointness.

A superinjective map from one curve complex to another is a map which
preserves disjointness and nondisjointness (superinjective maps are
easily seen to be simplicial and injective).  Theorem~\ref{intorelli} reduces to the following theorem.

\begin{thm}\label{superthm}
Let $g \geq 3$.  Every superinjective map $\Csep(S_g) \to \T(S_g)$ is
induced by an element of $\Mod(S_g)$.
\end{thm}

Let $S=S_g$, and let $\phi_\star : \Csep(S) \to \T(S)$ be a
superinjective map; in the original paper, $\phi_\star$ is induced by
an injective homomorphism $\phi:H \to \I(S)$, where $H$ is a finite
index subgroup of either $\I(S)$ or $\K(S)$.  The goal is to show
that $\phi_\star$ is induced by an element $f$ of $\Mod(S)$.  It follows that $f$ induces $\phi$, which gives Theorem~\ref{intorelli}.

A key idea for the argument in our original paper is that of a
sharing pair.  Suppose that $a$ and $b$ are separating curves
in $S$ which bound genus 1 subsurfaces $S_a$ and $S_b$ of $S$,
respectively.  We say that $a$ and $b$ form a \emph{sharing pair} for
the curve $\beta$ if $S-(S_a \cup S_b)$ is connected, and $S_a \cap
S_b$ is an annulus which contains the curve $\beta$.  Note that
$\beta$ is necessarily nonseparating.  The point is, if
a map $\phi_\star : \Csep(S) \to \Csep(S)$ preserves sharing pairs, then we can
extend $\phi_\star$ to a map $\C(S) \to \C(S)$: the nonseparating curve shared by $a$
and $b$ maps to the nonseparating curve shared by $\phi_\star(a)$ and
$\phi_\star(b)$.  

The basic outline of the proof of Theorem~\ref{superthm} is as follows.
\begin{enumerate}
\item The image of $\phi_\star$ lies in $\Csep(S)$.
\item The map $\phi_\star$ preserves topological types
  of curves, and it remembers when two curves are on the same side of
  a third curve.
\item The map $\phi_\star$ preserves sharing pairs.
\item $\phi_\star$ induces a well-defined superinjective map $\hat
  \phi_\star:\C(S) \to \C(S)$.
\item By a theorem of Irmak, $\hat \phi_\star$, hence $\phi_\star$, is induced by some $f
  \in \Mod(S)$ \cite{ei}.
\end{enumerate}
All of the arguments in the original paper are valid in the case of
genus 3 except the argument for Step 3.  In the remainder of this
addendum, we explain how to modify the proof of this step.

Step 3 above is Proposition 4.2 of the original paper, which is a
straightforward consequence of Steps 1 and 2 and of the following
lemma, which is Lemma 4.1 in the original paper. We restate the lemma
(with the case $g=3$ added) and explain how to modify the proof.

\begin{lem}
\label{nspchar}
Let $g \geq 3$, and let $a$ and $b$ be curves in $S=S_g$ which bound a
genus 1 subsurface of $S$.  Then $a$ and $b$ are a sharing pair if
and only if there exist separating curves $w$, $x$, $y$, and $z$ in
$S$ with the following properties.
\bl
\item $z$ bounds a genus 2 subsurface $S_z$ of $S$ 
\item $a$ and $b$ are in $S_z$ and intersect each other
\item $x$ and $y$ are disjoint
\item $w$ intersects $z$, but not $a$ and not $b$ 
\item $x$ intersects $a$ and $z$, but not $b$ 
\item $y$ intersects $b$ and $z$, but not $a$ 
\el  
\end{lem}

Note that we do not specify whether or not $w$
intersects $x$ or $y$.

One direction of the proof of Lemma~\ref{nspchar} works
as stated in the original paper for $g \geq 3$.  That is, if there
exist curves $w$, $x$, $y$, and $z$ with the given properties, then
$a$ and $b$ form a sharing pair.

\begin{figure}[h]
\psfrag{a}{$a$}
\psfrag{b}{$b$}
\psfrag{w}{$w$}
\psfrag{x}{$x$}
\psfrag{y}{$y$}
\psfrag{z}{$z$}
\centerline{\includegraphics[scale=.67]{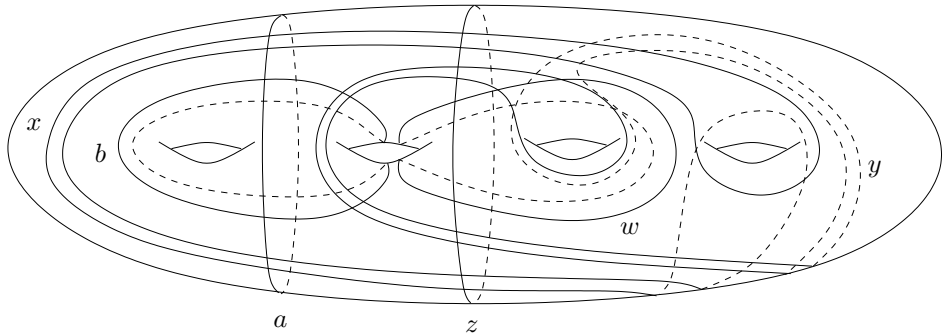}}
\caption{The original construction of the curves $w$, $x$, $y$, and $z$.}
\label{oldpic}
\end{figure}

It remains to show that, if $a$ and $b$ form a sharing pair in $S_g$
for $g \geq 3$, then we can find curves $w$, $x$, $y$, and $z$ which
satisfy the conditions of the lemma.  The idea from the original paper
is shown in Figure~\ref{oldpic}.  The reader will notice that this
construction is rather complicated, and does not give a useful
configuration in $S_3$ in any obvious way.

\begin{figure}[h]
\psfrag{a}{$a$}
\psfrag{b}{$b$}
\psfrag{w}{$w$}
\psfrag{x}{$x$}
\psfrag{y}{$y$}
\psfrag{z}{$z$}
\psfrag{...}{$\cdots$}
\centerline{\includegraphics[scale=.67]{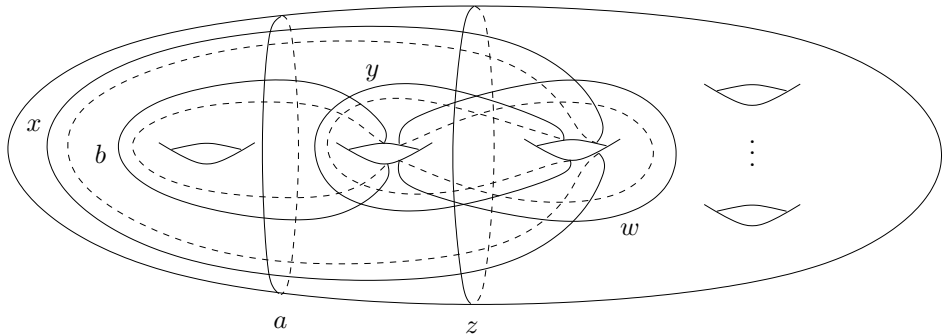}}
\caption{The new construction of the curves $w$, $x$, $y$, and $z$.}
\label{newpic}
\end{figure}

The new idea is to give a simpler configuration that works for every
genus $g \geq 3$.  This new configuration is shown in Figure~\ref{newpic}.

We remark that the configurations in Figures~\ref{oldpic}
and~\ref{newpic} are different in an essential way.  For instance,
consider the intersection of the genus 1 subsurface bounded by $y$ and
the genus 2 subsurface bounded by $z$; in Figure~\ref{oldpic}, this
intersection is a disk, whereas in Figure~\ref{newpic} it is an
annulus.  Perhaps more to the point, in Figure~\ref{oldpic}, each of
$w$, $x$, and $y$ bounds of a genus 1 subsurface, whereas in
Figure~\ref{newpic}, the curve $x$ does not bound a genus 1
subsurface when $g \geq 4$.

%%%
%%%
%%%

\bibliographystyle{plain}
\bibliography{kg2}

\end{document}